\documentclass[12pt]{amsart}

\theoremstyle{plain}
\newtheorem{thm}{Theorem}[section]
\newtheorem{theorem}[thm]{Theorem}

\newtheorem{lemma}[thm]{Lemma}

\newtheorem{proposition}[thm]{Proposition}
\theoremstyle{definition}
\newtheorem{remark}[thm]{Remark}

\newtheorem{definition}[thm]{Definition}

\newtheorem{assumption}[thm]{Assumption}

\newtheorem{example}[thm]{Example}

\newtheorem{question}[thm]{Question}

\numberwithin{equation}{section}

\newcommand{\sC}{{\mathcal C}}

\newcommand{\sF}{{\mathcal F}}

\newcommand{\sK}{{\mathcal K}}

\newcommand{\sO}{{\mathcal O}}
\newcommand{\sP}{{\mathcal P}}

\newcommand{\sT}{{\mathcal T}}

\newcommand{\sV}{{\mathcal V}}

\newcommand{\C}{{\mathbb C}}

\newcommand{\BP}{{\mathbb P}}

\title[Minimal rational curves]{Equivalence problem for minimal rational curves with isotrivial varieties of minimal rational tangents}
\author[Jun-Muk Hwang]{Jun-Muk Hwang${}^1$}
\address{Korea Institute for Advanced Study, Hoegiro 87, Seoul, 130-722, Korea} \email{jmhwang@kias.re.kr}
\thanks{${}^1$ supported by the Korea Research Foundation Grant
 (KRF-2006-341-C00004)}

\begin{document}

\maketitle

\begin{abstract}
We formulate the equivalence problem, in the sense of E. Cartan,
for families of minimal rational curves on uniruled projective
manifolds. An important invariant of this equivalence problem is
the variety of minimal rational tangents. We study the case when
varieties of minimal rational tangents at general points form an
isotrivial family.  The main question in this case is for which
projective variety $Z$, a family of minimal rational curves with
$Z$-isotrivial varieties of minimal rational tangents is locally
equivalent to the flat model. We show that this is the case when
$Z$ satisfies certain projective-geometric conditions, which hold
for a non-singular hypersurface of degree $\geq 4$.
\end{abstract}



\section{Introduction}

We will work over the complex numbers.    For a uniruled
projective manifold $X$,  an irreducible component $\sK$ of  the
space of rational curves on $X$ is  a {\it family of minimal
rational curves} on $X$  if the subvariety $\sK_x$ consisting of
members of $\sK$ through a general point $x \in X$ is projective
and non-empty. Minimal rational curves play an important role in
the geometry of uniruled projective manifolds (cf. [HM99], [Hw]).
We are interested in the following
 `equivalence problems' in the sense of E.  Cartan (c.f. \cite{Ca})
  for families of minimal
 rational curves.

\begin{question}\label{general Q} Let $X$ and $X'$ be two uniruled projective manifolds
with  families of minimal rational curves $\sK$  on $X$ and $\sK'$
on $X'$. Given two points $x \in X$ and $x' \in X'$, can we find
open neighborhoods $ x\in U \subset X$ and $x' \in U'\subset X'$
with a biholomorphic map $\varphi: U \to U'$ such that for each
member $C$ of $\sK$ (resp. $C'$ of $\sK'$) there exists a member
$C'$ of $\sK'$ (resp. $C$ of $\sK$) satisfying $$ \varphi(C \cap
U) = C' \cap U' \;\; \mbox{ (resp.
 } \varphi^{-1}(C' \cap U') = C \cap U)?$$ \end{question}

If such a biholomorphic map $\varphi$ exists, we will say that
$(X, \sK, x)$ is {\em equivalent } to $(X', \sK', x')$. One
motivation for studying this problem is the following theorem.

\begin{theorem}\label{CF} Let $X$ (resp. $X'$) be a Fano manifold with
second Betti number 1 and let $\sK$ (resp. $\sK'$) be a family of
minimal rational curves on $X$ (resp. $X'$).  Assume that $\dim
\sK = \dim \sK' \geq \dim X = \dim X'.$ Suppose for some $x \in X$
and $x' \in X'$, $(X, \sK, x)$ is equivalent to $(X', \sK', x')$.
Then the equivalence map $\varphi:U \to U'$ extends to a biregular
morphism from $X$ to $X'$ sending $x$ to $x'$.
\end{theorem}

Theorem \ref{CF} follows from the argument in \cite{HM01} although
it was not explicitly stated there. Theorem \ref{CF} and its
variations are useful  in proving two Fano manifolds of second
Betti number 1 are biregular (cf. [HM99], [Hw] ). Thus
 Question \ref{general Q} has interesting applications in
 algebraic geometry.

A natural approach to Question \ref{general Q} is  to find  local
properties of the family $\sK$ near $x$ which are invariant under
the equivalence, i.e., {\em local invariants} of the family. An
important invariant is provided by the variety of minimal rational
tangents. Recall that given a general point $x \in X$, the {\em
variety of minimal rational tangents} at $x$ is the subvariety
$\sC_x \subset \BP T_x(X)$ defined as the union of the tangent
directions of members of $\sK$ through $x$. A great advantage of
variety of minimal rational tangents $\sC_x$ is that it is
equipped with a projective embedding $\sC_x \subset \BP T_x(X)$
and consequently all projective geometric invariants of the
projective variety $\sC_x$ give rise to invariants of the
equivalence problem.

Throughout the paper we will consider only those $(X, \sK)$
   for which the following condition holds.
\begin{assumption}\label{assumption} $\dim X \geq 3$ and
$\sC_x$ at general point $x \in X$ is
   an irreducible non-singular variety and is not
   a linear subvariety in $\BP T_x(X)$. In particular, it has
   positive dimension.
\end{assumption}
What happens if $\sC_x$ is reducible is a very important and
difficult issue requiring ideas and methods different from those
considered below. One justification of making the assumption that
$\sC_x$ is irreducible is  that there is a large class of examples
satisfying it. As a matter of fact, all known examples with $\dim
\sC_x >0$ satisfy the irreducibility  assumption. The
non-singularity assumption is not really restrictive. It is
believed to be always true. Finally, the non-linearity assumption
is harmless.
 When $\sC_x$ is linear and irreducible, we can foliate $X$ by
projective spaces (e.g. [Ar, Theorem 3.1]) and the equivalence
problem becomes trivial.

 The main question in the
equivalence problem for minimal rational curves is to study to
what extent the equivalence is decided by the information of
varieties of minimal rational tangents.  More precisely, the main
question is the following.

\begin{question}\label{main Q}
Let $X$ and $X'$ be two  projective manifolds with families of
minimal rational curves $\sK$ on $X$ and $\sK'$ on $X'$ satisfying
Assumption \ref{assumption}. Let $x \in X$ and $x' \in X'$ be
general points in the sense of Assumption \ref{assumption}.
Suppose
 there exist open neighborhoods $  U \subset X$, $U' \subset X'$ and a commuting diagram
$$ \begin{array}{ccc} \BP T(U) & \stackrel{\psi}{\longrightarrow} & \BP T(U') \\\downarrow & & \downarrow \\
  U  & \stackrel{\varphi}{\longrightarrow} & U'\end{array}  $$
  where the vertical maps are natural
  projections and the
  horizontal maps are biholomorphisms satisfying $$\psi_x(\sC_x) =
  \sC_{\varphi(x)} \mbox{ for each } x \in U.$$ Is
  $(X, \sK, x)$  equivalent to $(X', \sK', x')$? \end{question}

We will see below that the answer is not affirmative in general. A
general result toward Question \ref{main Q} is provided by the
following result, which is just a restatement of Theorem 3.1.4 of
\cite{HM99}.

\begin{theorem}\label{HM99}
Let $X$ and $X'$ be two  projective manifolds with families of
minimal rational curves $\sK$ on $X$ and $\sK'$ on $X'$ satisfying
Assumption \ref{assumption}. Let $x \in X$ and $x' \in X'$ be
general points in the sense of Assumption \ref{assumption}.
Suppose
 there exist open neighborhoods $  U \subset X$, $U' \subset X'$ and a commuting diagram
$$ \begin{array}{ccc} \BP T(U) & \stackrel{\psi}{\longrightarrow} & \BP T(U') \\\downarrow & & \downarrow \\
  U  & \stackrel{\varphi}{\longrightarrow} & U'\end{array}  $$
  where the vertical maps are natural
  projections and the horizontal maps are biholomorphisms satisfying  $$\psi_x(\sC_x) =
  \sC_{\varphi(x)} \mbox{ for each } x \in U$$ and $\psi = d \varphi$, the derivative of $\varphi$. Then
  $(X, \sK, x)$ is  equivalent to $(X', \sK', x')$. \end{theorem}

In comparison to Question \ref{main Q},  the crucial additional
assumption in Theorem \ref{HM99} is that the holomorphic map
$\psi$ comes from the derivative of $\varphi$. In this sense, the
condition for Theorem \ref{HM99} is {\em differential-geometric}.
A central question is under what {\em algebraic-geometric}
conditions on the varieties of minimal rational tangents, we can
get this differential geometric condition. In this paper, we will
concentrate on the following special case.

\begin{definition} Let $Z \subset \BP^{n-1}$ be a fixed irreducible
non-singular non-linear projective variety. For an $n$-dimension
projective manifold $X$ and a family of minimal rational curves
$\sK$, we say that it has $Z$-isotrivial varieties of minimal
rational tangents, if for a general point $x \in X$, $\sC_x
\subset \BP T_x(X)$ is isomorphic to $Z \subset \BP^{n-1}$ as
projective varieties.\end{definition}

Note that for any $Z$, there exists $(X, \sK)$ with $Z$-isotrivial
varieties of minimal rational tangents:

\begin{example}\label{Ex} Let $Z \subset \BP^{n-1} \subset \BP^n$ be a non-singular
irreducible projective variety contained in a hyperplane. Let
$\psi: X_Z \to \BP^n$ be the blow-up of $\BP^n$ with center $Z$.
Let $\sK_Z$ be the family of curves which are proper transforms of
lines in $\BP^n$ intersecting $Z$. Then $\sK_Z$ is a family of
minimal rational curves on $X_Z$ with $Z$-isotrivial varieties of
minimal rational tangents. In fact, $X_Z$ is quasi-homogeneous
with an open orbit containing $\psi^{-1}(\BP^n \setminus
\BP^{n-1})$. \end{example}

Now we can formulate the following special case of Question
\ref{main Q}.

\begin{question}\label{iso Q}
Let $Z \subset \BP^{n-1} $ be an irreducible non-singular
non-linear variety. Let $X$ be an $n$-dimensional projective
manifold  and let $\sK$ be a family of minimal rational curves on
$X$ with $Z$-isotrivial varieties of minimal rational tangents. Is
$(X, \sK, x)$ for a general $x \in X$  equivalent to that of
Example \ref{Ex}?
\end{question}

The answer is not always affirmative:

\begin{example}\label{counter} Let $W$ be a $2 \ell$-dimensional
complex vector space with a symplectic form. Fix an integer $k$,
$1 < k < \ell$ and let $S$ be the variety of all $k$-dimensional
isotropic subspaces of $W$. $S$ is a uniruled homogeneous
projective manifold. There is a unique family  $\sK$ of minimal
rational curves, just the set of all lines on $S$ under the
Pl\"ucker embedding. The varieties of minimal rational tangents
are $Z$-isotrivial where $Z$ is the projectivization of the vector
bundle $\sO(-1)^{2 \ell -2k} \oplus \sO(-2)$ on $\BP^{k-1}$
embedded by the dual tautological bundle of the projective bundle
(cf. Proposition 3.2.1 of \cite{HM05}). Let us denote it  by $Z
\subset \BP V$. There is a distinguished hypersurface $R \subset
Z$ corresponding to $\BP \sO(-1)^{2\ell -2k}$. Let $D$ be the
linear span of $R$ in $V$. This $D$ defines a distribution on $S$
which is not integrable (cf. Section 4 of \cite{HM05}). However,
the corresponding distribution on $X_Z$ of Example \ref{Ex} is
integrable. Thus $(S, \sK, x)$ cannot be equivalent to $(X_Z,
\sK_Z, y)$ at general points $x, y$.
\end{example}

Thus the correct formulation of Question \ref{iso Q} is to ask for
which $Z$ the answer to Question \ref{iso Q} is affirmative.
 Up to now the
only result in this line is the following result of Ngaiming Mok
in \cite{Mo}:

\begin{theorem}\label{Mok}
Let $S$ be an $n$-dimensional irreducible Hermitian symmetric
space of compact type with a base point $o \in S$. If the
projective  variety $Z \subset \BP^{n-1}$ is isomorphic to  $
\sC_o \subset \BP T_o(S)$ for the family of minimal rational
curves on $S$, then Question \ref{iso Q} has an affirmative
answer. \end{theorem}

For example when $S$ is the $n$-dimensional quadric hypersurface,
$Z \subset \BP^{n-1}$ is just an $(n-2)$-dimensional non-singular
quadric hypersurface. Then $\sC_x \subset \BP T_x(X)$ in Question
\ref{iso Q} defines a conformal structure at general points of
$X$. In this case, Theorem \ref{Mok} says that this conformal
structure is flat. In general, for each $S$, we can interpret the
condition of Question \ref{iso Q} as  a certain G-structure at
general points of $X$ and Theorem \ref{Mok} says that this
G-structure is flat.

It is worth recalling Mok's strategy for the proof of Theorem
\ref{Mok}. The main point is to show that this G-structure which
is defined at general points of $X$ can be extended to a
G-structure in a {\em neighborhood } of a general minimal rational
curve. Once this extension is obtained, one can deduce the
flatness by applying \cite{HM97} which shows the vanishing of the
curvature tensor from {\em global} information of the tangent
bundle of $X$ on the minimal rational curve.

The projective variety $Z \subset \BP^{n-1}$ treated in Theorem
\ref{Mok} is a homogeneous variety with reductive automorphism
group.  Our main result concerns the opposite case when the
automorphisms of $Z \subset \BP^{n-1}$ is 0-dimensional:

\begin{theorem}\label{Main T}  Assume that  $Z \subset \BP V$ satisfies the following
conditions.

(1) $Z$ is non-singular, non-degenerate and linearly normal.

(2) The variety of tangent lines to $Z$, defined as a subvariety
of ${\bf Gr}( 2, V) \subset \BP (\wedge^2 V)$, is non-degenerate
in $\BP (\wedge^2 V)$.

 (3) $H^0(Z, T(Z)
\otimes \sO (1)) = H^0(Z, {\rm ad}(T(Z)) \otimes \sO (1)) =0$
where ${\rm ad}(T(Z))$ denotes the bundle of traceless
endomorphisms of the tangent bundle of $Z$.

Then for  a uniruled manifold $X$ and a family $\sK$ of minimal
rational curves with $Z$-isotrivial varieties of minimal rational
tangents, $(X, \sK, x)$ at a general point $x$ is equivalent to
that of Example \ref{Ex}.
\end{theorem}

Note that the non-degeneracy and $H^0(Z, T(Z) \otimes \sO (1)) =0$
imply that the projective automorphism group of $(Z \subset \BP
V)$ is $0$-dimensional. This means that we have a G-structure at
general points of $X$ with the group G isomorphic to the scalar
multiplication group $\C^*$. The essential difference from Theorem
\ref{Mok} is the following: the G-structure {\em cannot} be
extended to a neighborhood of a minimal rational curves. One way
to see this is by directly checking it in the case of Example
\ref{Ex}. There is a more conceptual way to see this as follows.
Suppose it is possible to extend the G-structure to a neighborhood
$U$ of a general minimal rational curve. For simplicity, let us
assume that $U$ is simply connected and $Z \subset \BP T(U)$ has
no non-trivial automorphism. Then in $\BP T(U)$ we have a
submanifold $\sC \subset \BP T(U)$ with each fiber $\sC_x \subset
\BP T_x(U)$ isomorphic to $Z \subset \BP^{n-1}$. Since the
automorphism group of $Z \subset \BP^{n-1}$ is trivial, we get a
unique trivialization of the projective bundle $\BP T(U)$. But on
a general minimal rational curve, $T(U)$ splits into $\sO(2)
\oplus \sO(1)^{p} \oplus \sO^{n-1-p}$ for some $p>0$, a
contradiction. Thus in the case where the automorphism group of
$Z$ is 0-dimensional, we cannot use Mok's approach. In the setting
of Theorem \ref{Main T}, the flatness of the G-structure, or the
vanishing of the corresponding `curvature tensors', must be proved
only at general points. In other words, it must come from
information on the geometry of minimal rational curves in a
neighborhood of a general point. The crucial point of the proof of
Theorem \ref{Main T} lies in the use of  such local information
 to prove  flatness.

The variety of minimal rational tangents $Z \subset \BP V$ in
Example \ref{counter}, satisfy the conditions (1) and (2) in
Theorem \ref{Main T}. Thus some additional conditions like (3) are
necessary. However we expect that the condition  (3) can be
weakened, hopefully, to $H^0(Z, T(Z)) = H^0(Z, {\rm ad}(T(Z)))
=0$.

Examples of $Z \subset \BP V$ satisfying the conditions (1)-(3) of
Theorem \ref{Main T} are provided by non-singular hypersrufaces of
degree $\geq 4$ in $\BP V$, $ \dim V \geq 3$. In fact, (1) is
standard for hypersurfaces and (2) follows from Proposition 2.6 of
\cite{Hw}. (3) can be checked from Theorem 4 (iii) of \cite{BR}.
It is likely that the three conditions for $Z$ also hold for a
large class of complete intersections $Z \subset \BP V$.

\section{Coframes on a manifold and induced coframes on the tangent bundle}

This section and the next section  are  concerned with the basic
coframe formulation of Cartan's approach to equivalence problems.
Its content must be well-known to differential geometers and
essentially covered by [Ca]. However the special case we need had
not been explicitly worked out. Also we expect some of the readers
have background in algebraic geometry. So we will give a
self-contained presentation with more coordinate-free notation.
  All functions, differential forms and tensors
considered here are assumed to be holomorphic.

Let us start by recalling the notion of differential forms with
values in a vector space.  Section V.6 of \cite{St} is a good
reference for this notion.  Given a vector space $V$, a $V$-valued
differential $k$-form on $M$ is just a holomorphic section $\omega$
of the vector bundle $\Omega^k_M \otimes V_M$ where $V_M$ denotes
the trivial bundle on $M$ with each fiber $V$. The exterior
derivative $d \omega$ is defined as a $V$-valued $(k+1)$-form and
satisfies the usual properties of the exterior derivative including
$d(d \omega)=0$. For a $V_1$-valued $k_1$-form $\omega_1$ and
$V_2$-valued $k_2$-form $\omega_2$, their exterior product $\omega_1
\wedge \omega_2$ is defined as $(V_1 \otimes V_2)$-valued $(k_1 +
k_2)$-form. Let $W$ be another vector space and $\rho$ be a ${\rm
Hom}(V, W)$-valued function on $M$. Given any $V$-valued $k$-form
$\omega$, $\rho_{\sharp} \omega$ is the $W$-valued $k$-form defined
by the composition of $\omega$ with $\rho$. For another vector space
$W'$ and ${\rm Hom}(W, W')$-valued function $\eta$ on $M$, we have
$$ (\eta \circ \rho)_{\sharp} \omega = \eta_{\sharp} ( \rho_{\sharp}
\omega).$$ Similarly, we have the notion of $V$-valued vector field
$D$, as a holomorphic section of $T(M) \otimes V_M$. For a ${\rm
Hom}(V, W)$-valued function $\rho$ on $M$, $\rho_{\sharp} D$ is a
$W$-valued vector field. For a $V$-valued vector field $D$ and a
$W$-valued 1-form $\omega$, $D \rfloor \omega$ denotes the $(V
\otimes W)$-valued function obtained by the natural pairing.

\begin{definition}
Let $M$ be a complex manifold of dimension $n$. Fix a vector space
$V$ of dimension $n$. A $V$-valued 1-form $\omega$ on $M$ is called
a {\em coframe} if for each $x \in M$, the homomorphism $\omega_x :
T_x(X) \rightarrow V$ is an isomorphism. \end{definition}

 The
following is immediate from a point-wise consideration.

\begin{lemma}\label{>2} Let $\omega$ be a coframe on a manifold $M$ of dimension $\geq 3$.
\begin{enumerate} \item[1.]   Let $\xi$ and $\xi'$ be two
$V^*$-valued functions with $\xi_{\sharp} \omega = \xi'_{\sharp}
\omega.$ Then $\xi = \xi'$. \item[2.] Let $\xi$ and $\xi'$ be two
  $W$-valued $k$-forms, $k =
1,2$, for some finite-dimensional vector space $W$. Suppose that
$\xi \wedge \omega= \xi' \wedge  \omega.$ Then $\xi = \xi'$.
\item[3.] Let $\xi$ and $\xi'$ be two
  $W$-valued functions,  for some finite-dimensional vector space $W$.
Suppose that $\xi_{\sharp}( \omega \wedge \omega) =
\xi'_{\sharp}(\omega \wedge \omega)$. Then $\xi = \xi'$.
\end{enumerate}
\end{lemma}

\begin{definition} Given a coframe $\omega$, the {\em dual frame} of $\omega$ is
a  $V^*$-valued vector field $D_{\omega}$ on $M$, whose value at $x
\in M$ is given by the inverse of the isomorphism $\omega_x: T_x(X)
\rightarrow V$. In other words, the $(V^* \otimes V)$-valued
function $D_{\omega} \rfloor \omega$ has constant value $ {\rm Id}_V
\in V^* \otimes V.$
\end{definition}

  The following can be checked easily.

\begin{lemma}\label{d-lemma} For a coframe $\omega$ and its dual frame
$D_{\omega}$, the following holds.
\begin{enumerate} \item[1.]  For any function $f$, we have
$$df = (D_{\omega} f)_{\sharp} \omega$$ where $D_{\omega}f$ is the
$V^*$-valued function  obtained by differentiating $f$ by
$D_{\omega}$. \item[2.] For any vector field $v$ on $M$, the
$V$-valued function $v \rfloor \omega$ satisfies
$$(v \rfloor \omega)_{\sharp} D_{\omega} = v.$$  \end{enumerate}\end{lemma}

\begin{definition} A coframe $\omega$ is {\em closed} if $d \omega =0$.
 Given a coframe $\omega$, there exists a ${\rm Hom}( V
\otimes V, V)$-valued function $\sigma^{\omega}$ on $M$, called
the {\em structure function } of $\omega$  such that
$$ d \omega = \sigma^{\omega}_{\sharp}(\omega \wedge \omega).$$
In fact, $\sigma^{\omega}$ takes values in ${\rm Hom}(\wedge^2 V,
V) \subset {\rm Hom}(V \otimes V, V)$ from the anti-symmetry in $d
\omega$. By the canonical isomorphism ${\rm Hom}( V \otimes V, V)
= {\rm Hom}(V^*, V^* \otimes V^*)$, we can view the structure
function as a ${\rm Hom}(V^*, V^* \otimes V^*)$-valued function,
in which case we will denote it by $\delta^{\omega}$. Then $$ [
D_{\omega}, D_{\omega}] = \delta^{\omega}_{\sharp} D_{\omega}.$$
\end{definition}

A coframe $\omega$ on a manifold $M$ induces a coframe $\Omega$ on
the manifold $T(M)$ in the following way.

\begin{definition} View a coframe $\omega$ on $M$ as a $V$-valued function $\mu$ on $T(M)$.
Let $\pi: T(M) \to M$ be the projection and let $\theta:= \pi^*
\omega$ be the $V$-valued 1-form on $T(M)$ obtained by pulling back
$\omega$. Then the pair
$$ \Omega : = (\theta, \lambda:= d \mu)$$ is a $(V \oplus
V)$-valued 1-form on $T(M)$, which is in fact a coframe on $T(M)$.
This is the {\em induced coframe} on $T(M)$. \end{definition}

\begin{proposition}\label{Omega-str}   To avoid confusion, we will rename $V \oplus V$ as $V_1
\oplus V_2$. Then the structure function of the coframe $\Omega$,
which is a ${\rm Hom}(\wedge^2(V_1 \oplus V_2), V_1 \oplus
V_2)$-valued function on $T(M)$ takes values in ${\rm Hom} (\wedge^2
V_1, V_1)$. In fact,$$\sigma^{\Omega} = \pi^* \sigma^{\omega}.$$
\end{proposition}

\begin{proof}
This is immediate from
$$ d \Omega = (d \theta, d \lambda) = (\pi^* d \omega, 0).$$ \end{proof}

The following is straight forward from the definitions and
Proposition \ref{Omega-str}.

\begin{proposition}\label{dual}  The dual frame of the coframe $\Omega$  is
given by
$$D_{\Omega} = ( D_{\theta}, D_{\lambda})$$ where $D_{\theta}$ and
$D_{\lambda}$ are $V^*$-valued vector fields satisfying
\begin{equation}\label{dual-eq} D_{\theta} \rfloor \theta = D_{\lambda} \rfloor \lambda =
D_{\lambda}\mu = {\rm Id}_V \mbox{ and } D_{\theta}\mu= D_{\theta}
\rfloor \lambda = D_{\lambda} \rfloor \theta = 0.
\end{equation} They have the following properties.
\begin{enumerate}
\item[(a)] Under the projection $d \pi: T(T(M)) \to T(M)$,
 $d \pi( D_{\theta}) = D_{\omega}.$
 \item[(b)] Any vector field $\tilde{v}$ on an open subset $U \subset T(M)$ tangent to the fibers
 of $\pi$ is of the form $\tilde{v} = h_{\sharp} D_{\lambda} \mbox{ for some $V$-valued
  function $h$ on $U$.}$
  \item[(c)]
$[D_{\theta}, D_{\lambda}] = [D_{\lambda}, D_{\lambda}] = 0, \;\;
[D_{\theta}, D_{\theta}] = (\pi^*\delta^{\omega})_{\sharp}
D_{\theta}.$
\end{enumerate} \end{proposition}

\begin{definition} For a coframe $\omega$ on a manifold $M$,
 the vector field on $T(M)$ $$\gamma := \mu_{\sharp}
D_{\theta}$$ is called the {\em geodesic flow} of the coframe
$\omega$.
\end{definition}

\begin{proposition}\label{gamma} Given a coframe $\omega$ on $M$,
its geodesic flow $\gamma$ on $T(M)$ has the following properties.
\begin{enumerate}
\item[(a)] $[D_{\lambda}, \gamma] = D_{\theta}.$ \item[(b)] For a
point $v \in T(M)$, let $\gamma_v \in T_{v}(T(M))$ be the value of
$\gamma$ at $v$. Then
$$  d \pi_v (\gamma_v) = v \mbox{ for any vector } v \in T(M).$$
   \end{enumerate} \end{proposition}

\begin{proof}
From Proposition \ref{dual} (c) and \eqref{dual-eq},
$$[D_{\lambda}, \gamma] = [ D_{\lambda}, \mu_{\sharp} D_{\theta}] =
(D_{\lambda}\mu)_{\sharp} D_{\theta} = D_{\theta},$$ proving (a).

For (b), let $\pi(v) = x \in M$.  From $$ (\mu_{\sharp}
D_{\theta})_v = \mu(v)_{\sharp} (D_{\theta})_v =
\omega(v)_{\sharp} (D_{\theta})_v$$ and $d \pi_v(D_{\theta})_v =
(D_{\omega})_x$,
$$ d \pi_v (\gamma_v) = d \pi_v ((\mu_{\sharp} D_{\theta})_v)
=  \omega(v)_{\sharp} (D_{\omega})_x.$$ The latter must be $v$ by
Lemma \ref{d-lemma}.
\end{proof}

\section{ Conformally closed coframes }

Let $\omega$ be a coframe on a manifold $M$.
 For any function $h$ on $M$, $(D_{\omega} h)_{\sharp}\omega$ is a
 closed 1-form by Lemma \ref{d-lemma} (a). The following converse is
just a restatement of the Poincar\'e lemma.

\begin{lemma}\label{Poincare}
Given a coframe $\omega$ and
 a $V^*$-valued function $\xi$ on $M$, suppose that $\xi_{\sharp} \omega$ is a
closed 1-form. Then for any point $x \in M$, there exists a
neighborhood $x\in U\subset M$ and a function $h$ on $U$ with $\xi
= D_{\omega} h$.
\end{lemma}

\begin{definition} A coframe $\omega$ is {\em
conformally closed} if for any point $x \in M$, there exist a
neighborhood $x \in U \subset M$  and a non-vanishing  function $f
$ on $U$ such that $f \omega $ is closed on $U$.
\end{definition}

We will need the following elementary fact from linear algebra.
One can check it either by an explicit computation or from the
fact that ${\rm Hom}(\wedge^2 V, V)$ decomposes into two
irreducible factors as a ${\rm GL}(V)$-module.

\begin{proposition}\label{Xi} For a vector space $V$, let $\Xi_V \subset {\rm
Hom}(\wedge^2 V, V)$ be the subspace defined by $$ \Xi_V := \{
\sigma: V\otimes V \to V, \sigma(u, v) = - \sigma(v, u) \in \C u +
\C v \mbox{ for any } u, v \in V \}.$$ Define the contraction
homomorphism $\iota: V^* \to {\rm Hom}(V \otimes V, V)$ such that
for a $V^*$-valued function $\eta$, $$ \iota(\eta)_{\sharp}
(\omega \wedge \omega) = (\eta_{\sharp} \omega) \wedge \omega.$$
Then $\iota$ is injective and its image is $\Xi_V$.
\end{proposition}

 The following theorem characterizes
conformally closed coframes in terms of their structure functions.

\begin{theorem}\label{XiV} A coframe $\omega$ on a  manifold of dimension
$\geq 3$ is conformally closed if and
only if its structure function $\sigma_{\omega}$ takes  values in
$\Xi_V \subset {\rm Hom}(\wedge^2 V, V)$. \end{theorem}

\begin{proof} Suppose that $f \omega$ is closed, i.e., $$0= d (f \omega) =
df \wedge \omega + f d \omega.$$ Then
$$(\sigma_\omega)_{\sharp}(\omega \wedge \omega) = d \omega =
-\frac{df}{f}  \wedge \omega = -d ( \log f) \wedge \omega = -
((D_{\omega}(\log f))_{\sharp} \omega) \wedge \omega.$$ This is
equal to $\iota(-D_{\omega}(\log f))_{\sharp} (\omega \wedge
\omega).$ By Lemma \ref{>2}, we conclude that $$ \sigma^{\omega} =
\iota(-D_{\omega}(\log f)).$$ Thus $\sigma^{\omega}$ takes values
in $\Xi_V$.

Conversely, assume that $\sigma^{\omega}$ takes values in $\Xi_V$,
i.e. there exists some $V^*$-valued function $\xi$ on $M$  such
that
$$d \omega = \iota(\xi)_{\sharp} (\omega \wedge \omega)=
(\xi_{\sharp} \omega ) \wedge \omega.$$ Taking $d$-derivative on
both sides, we get
$$0 = d d \omega = d (\xi_{\sharp}\omega)  \wedge \omega  + (\xi_{\sharp} \omega)
\wedge  d \omega.$$ Note that
$$ (\xi_{\sharp} \omega)
\wedge  d \omega = (\xi_{\sharp} \omega) \wedge
\iota(\xi)_{\sharp} (\omega \wedge \omega) = (\xi_{\sharp} \omega)
\wedge (\xi_{\sharp} \omega) \wedge \omega.$$ Since $\xi_{\sharp}
\omega$ is a (scalar-valued) 1-form on $M$, the right hand side
must be 0. It follows that $d (\xi_{\sharp}\omega)  \wedge \omega
=0$. From Lemma \ref{>2}, this implies that $\xi_{\sharp}\omega$
is closed. Then by Lemma \ref{Poincare}, there exist a
neighborhood $U$ and a function $h$ on $M$ such that $\xi =
D_{\omega} h$ and, via Lemma \ref{d-lemma},
$$ d \omega = (\xi_{\sharp}\omega) \wedge \omega
= ((D_{\omega} h)_{\sharp} \omega) \wedge \omega = d h \wedge
\omega.$$ Write $ h = - \log f$. Then
$$ d (f \omega) = df \wedge \omega + f d \omega = df \wedge \omega + f dh \wedge \omega = df \wedge \omega
+ f (-\frac{df}{f}) \wedge \omega = 0.$$ Thus $\omega$ is
conformally closed.
\end{proof}

\begin{remark} Theorem \ref{XiV} is false when the dimension of
the manifold is 2. In fact, when $\dim V =2$, $\Xi_V = {\rm
Hom}(\wedge^2V, V)$ and the condition in Theorem \ref{XiV} is
empty. \end{remark}

\section{Coframes adapted to an isotrivial cone structure  }

\begin{definition} A {\em cone structure } on a complex manifold
$M$ is a submanifold  $\sC \subset \BP T(M)$ such that the
projection $\varpi: \sC \to M$ is a smooth morphism with connected
fibers. For each point $x \in M$, the fiber $\varpi^{-1}(x)$ will
be denoted by $\sC_x. $
\end{definition}

\begin{definition}  Let $V$ be an $n$-dimensional vector space and let $Z \subset
\BP V$ be a fixed projective subvariety. A cone structure $\sC
\subset \BP T(M)$ on an $n$-dimensional manifold $M$ is said to be
$Z$-{\em isotrivial} if for each $ x\in M$, the inclusion $(\sC_x
\subset \BP T_x(M))$ is isomorphic to $(Z \subset \BP V)$ up to
projective transformations. \end{definition}

\begin{definition}
Given $Z \subset \BP V$ and a $Z$-isotrivial cone structure on
$M$, a coframe $\omega$ on $M$ is said to be {\it adapted} to the
cone structure if for each $x\in M$, the isomorphism $\omega_x:
T_x(M) \to V$ sends $\sC_x \subset \BP T_x(M)$ to $Z \subset \BP
V$. Given any $Z$-isotrivial cone structure on a manifold $M$, an
adapted coframe exists if we shrink $M$.
\end{definition}

\begin{definition}\label{flat}
A $Z$-isotrivial cone structure $\sC \subset T(M)$ is {\em locally
flat} if for any point $x \in M$ there exist a neighborhood $x \in
U \subset M$ and  a biholomorphic map $\zeta: U \to V$ such that
$$ \zeta_*(\sC|_U) = \zeta(U) \times Z \subset V \times \BP V$$
where $  \zeta_*: \BP T(M) \to \BP T(V) = V \times \BP V$ is the
differential of $\zeta$.  \end{definition}

\begin{proposition} A $Z$-isotrivial cone structure is locally
flat if and only if it has a conformally closed adapted frame.
\end{proposition}

\begin{proof} The `only if' part is obvious: the differential of the
map $\zeta$ in Definition \ref{flat} defines a closed adapted
coframe. For the `if' part, since the question is local, we may
assume that there  is a  closed adapted coframe. Then by
Poincar\'e lemma, we can integrate it in a neighborhood to get a
biholomorphic map $\zeta:U \to V$ satisfying the required
condition.
\end{proof}

We  collect a few properties of adapted coframes.

\begin{proposition}\label{gamma-tan}
Given a $Z$-isotrivial cone structure $\sC \subset \BP T(M)$, let
$\hat{\sC} \subset T(M)$ be the cone over $\sC$. For an adapted
coframe $\omega$, the  cone $\hat{\sC} \subset T(M)$ is preserved
by the geodesic flow $\gamma$ of $\omega$, i.e., at any point $u
\in \hat{\sC}\setminus 0$, $\gamma_u \in T_u(\hat{\sC})$.
\end{proposition}

\begin{proof}
Let $\mu$ be the $V$-valued function on $T(M)$ defined by
$\omega$. Let $I_Z \subset {\rm Sym}^{\bullet} V^*$ be the ideal
defining the projective variety $Z$.  Since $\omega$ is adapted,
$\sC \subset T(M)$ is defined as the zero locus of the collection
of  functions
$$\{ g\circ \mu: T(M) \stackrel{\mu}{\longrightarrow} V
\stackrel{g}{\longrightarrow} \C, \; g \in I_Z \}.$$ Since
$D_{\theta} \mu = 0$, we see that $$\gamma (g \circ \mu)  =
\mu_{\sharp} D_{\theta} (g \circ \mu) =0.$$ Thus $\gamma$ is
tangent to $\hat{\sC}$.
\end{proof}

\begin{proposition} \label{isoV}
Let $\omega$ be a coframe adapted to a $Z$-isotrivial cone
structure $\sC \subset \BP T(M)$. Let $u$ be a non-zero point of
the affine cone $\hat{\sC} \subset T(M)$ and $ v \in
T_{u}(\hat{\sC}_x)$ with $ x = \pi(u) \in M,$. Then there exists a
local $V$-valued function $g$ in a neighborhood of $u$ in $T(M)$
satisfying $D_{\theta}g = 0$ such that the vector field
$g_{\sharp} D_{\lambda}$ is tangent to the fibers of $\hat{\sC}
\to M$ and $(g_{\sharp} D_{\lambda})_u = v.$ \end{proposition}

\begin{proof} Let $\bar{u} = \mu(u) \in \hat{Z} \subset V$ and
$$\mu_*: T_u(T(M)) \to T_{\bar{u}}(V)$$ be the differential of the function $\mu:
T(M) \to V$. Choose a germ of vector field $\vec{v}$ at $\bar{u}
\in V$ such that $\vec{v}$ is tangent to $\hat{Z}$ and
$\vec{v}_{\bar{u}} = \mu_* (v)$. Under the canonical
trivialization of $T(V)= V \times V$, the vector field $\vec{v}$
defines a $V$-valued function $\bar{g}$ in a neighborhood of
$\bar{u}$. The $V$-valued function $g:= \mu^* \bar{g}$ in a
neighborhood of $u$ satisfies $D_{\theta} g = 0$ from
\eqref{dual-eq} in Proposition \ref{dual}. The vector field
$g_{\sharp} D_{\lambda}$ is tangent to the fibers of $\hat{\sC}
\to M$ and its value at $u$ is $v$, from the choice of $\bar{g}$.
\end{proof}

\section{Characteristic connection}

\begin{definition} Given a cone structure $\sC \subset \BP T(M)$,
denote by $\varpi: \sC \to M$ the natural projection.  Denote by
$\sV \subset T(\sC)$ the relative tangent bundle of the projection
$\varpi$ and by   $\sT \subset T(\sC)$ the tautological bundle
whose fiber at $\alpha \in \sC_x$ is $d
\varpi_{\alpha}^{-1}(\hat{\alpha}) $ where $\hat{\alpha} \subset
T_x(M)$ is the 1-dimensional subspace corresponding to $\alpha \in
\BP T_x(M)$.
 A line subbundle $\sF \subset T(\sC)$, with locally free quotient
$ T(\sC)/ \sF$, is called a {\em conic connection}  if $\sF
\subset \sT$ and $\sF \cap \sV =0$, i.e., it splits the exact
sequence \begin{equation} \label{split} 0 \longrightarrow \sV
\longrightarrow \sT \longrightarrow \sT/\sV \cong \sO(-1)
\longrightarrow 0\end{equation} where $\sO(1)$ denotes the
relative hyperplane bundle on $\BP T(M)$.
\end{definition}

\begin{proposition}\label{unique} Let $\sC \subset \BP T(M)$ be a cone structure.
 Suppose that $H^0(\sC_x, \sV \otimes \sO(1)) =0$ for some $x \in M$. Then a conic connection is unique if it exists.
\end{proposition}

\begin{proof} This follows from the fact that the set of splittings
of \eqref{split} is $H^0(\sC, \sV \otimes \sO(1))$.
\end{proof}

Given a vector space $V$ and a non-singular projective variety $Z
\subset \BP V$, the cone over $Z$ will be denoted by $\hat{Z}
\subset V$. For a point $\alpha \in Z$,  the affine tangent space
of $Z$ at $\alpha$ is $$\hat{T}_{\alpha}(Z):= T_u(\hat{Z}) \subset
V \mbox{ a non-zero vector } u \in \hat{\alpha}.$$ This is
independent of the choice of $u$.

 Let $ \sC \subset \BP T(M)$ be a
cone structure. Denote by $\sP \subset T(\sC)$ the subbundle whose
fiber at $\alpha \in \sC_x$ is $d
\varpi_{\alpha}^{-1}(\hat{T}_{\alpha}(\sC_x))$ where
$\hat{T}_{\alpha}(\sC_x) \subset T_x(M)$ is the affine tangent
space of the projective subvariety $\sC_x \subset \BP T_x(M)$ at
$\alpha \in \sC_x$. The following is proved in Proposition 1 of
\cite{HM04}.

\begin{proposition}\label{sP} Given a conic connection $\sF$, regarding the subbundles  of
$T(\sC)$ as sheaves of vector fields on $\sC$, we have $\sP =
[\sF, \sV].$
\end{proposition}

\begin{proposition}\label{Gamma}
Given a $Z$-isotrivial cone structure $\sC \subset \BP T(M)$ and
an adapted coframe $\omega$, the geodesic flow $\gamma$ is tangent
to the cone $\hat{\sC} \subset T(M)$ by Proposition
\ref{gamma-tan}. Denote by $\Gamma \subset T(\sC)$ the line
subbundle spanned by the image of $\gamma$.  Then $\Gamma$ is a
conic connection on the cone structure $\sC \subset \BP T(M)$.
\end{proposition}

\begin{proof} It suffices to show $\Gamma \subset \sT$. This is immediate from
Proposition \ref{gamma} (b). \end{proof}

Let us introduce a distinguished class of conic connections.

\begin{definition} A conic connection $\sF \subset T(\sC)$ is a {\em
characteristic connection} if  for any local section $v$ of $\sP$
and any local section $w$ of $\sF$, both regarded as local vector
fields on the manifold $\sC$, the Lie bracket $[v, w]$ is a local
section of $\sP$ again. \end{definition}

When  $\sC = \BP T(M)$, any conic connection is a characteristic
connection.  On the other hand, when  $\sC \neq \BP T(M)$,  a
characteristic connection is unique if it exists (Theorem 3.1.4 of
\cite{HM99}).

\begin{proposition}\label{XiZ} Given a non-singular projective variety
 $Z \subset \BP V$, define the subspace $\Xi_Z \subset {\rm
Hom}(\wedge^2 V, V)$ by $$\Xi_Z := \{ \sigma: V \otimes V \to V,
\; \sigma(u, v) = - \sigma(v, u), \sigma(u, v) \in
\hat{T}_{\alpha}(Z) \mbox{ if } \alpha \in Z,  u \in \hat{\alpha}
\mbox{ and } v \in \hat{T}_{\alpha}(Z) \}.$$ Let $\sC \subset \BP
T(M)$ be a $Z$-isotrivial cone structure and $\omega$ be an
adapted coframe. Suppose the conic connection $\Gamma$ on $\sC$
induced by $\omega$ in the sense of Proposition \ref{Gamma} is a
characteristic connection. Then the structure function
$\sigma^{\omega}$ of $\omega$ takes values in $\Xi_Z$.
\end{proposition}

\begin{proof}
By Proposition \ref{sP}, if $\sF$ is a characteristic connection,
then given any local section $\tilde{v}$ of $\sV$ and a local
section $\tilde{w}$ of $\sF$, $[[\tilde{v}, \tilde{w}],
\tilde{w}]$ is a local section of $\sP$. We can pull-back this to
the affine cone $\hat{\sC} \subset T(M)$ by the projection
$\hat{\sC} \setminus 0 \to \sC$. Let us denote the pull-back
distributions of $\sV$ and $\sP$ by $\hat{\sV}$ and $\hat{\sP}$.
By Proposition  \ref{isoV}, for any vector $v \in
\hat{T}_{\alpha}(\hat{\sC_x})$, there exists a section $\tilde{v}$
of the relative tangent bundle $\hat{\sV}$ of $\hat{\sC} \to M$
such that
$$ \tilde{v}_u = v, \; \tilde{v} = g_{\sharp} D_{\lambda} \mbox{
for some $V$-valued function } g \mbox{ satisfying } D_{\theta} g
= 0.$$  From $\gamma(g) = 0$ and Proposition \ref{gamma} (a),
$$[\tilde{v}, \gamma] = [g_{\sharp} D_{\lambda}, \gamma] =
g_{\sharp} [D_{\lambda}, \gamma] - \gamma(g)_{\sharp} D_{\lambda}
= g_{\sharp} D_{\theta}.$$
$$[[\tilde{v}, \gamma], \gamma] = [g_{\sharp} D_{\theta}, \gamma]
= g_{\sharp} [D_{\theta}, \gamma] = g_{\sharp} [D_{\theta},
\mu_{\sharp} D_{\theta}] = g_{\sharp} (\mu_{\sharp} [ D_{\theta},
D_{\theta}]) + g_{\sharp} (D_{\theta} \mu) D_{\theta}.$$ Since
$D_{\theta} \mu = 0$, we see that
$$d \pi_u ([[\tilde{v}, \gamma],
\gamma]_u) = g(u)_{\sharp}( \mu(u)_{\sharp} [D_{\theta},
D_{\theta}]_u)= \sigma^{\omega}(v, u).$$ Since $\Gamma$ is a
characteristic connection, the latter must be in
$\hat{T}_{\alpha}(Z)$ by Proposition \ref{sP}. Thus
$\sigma^{\omega} \in \Xi_Z$.
\end{proof}

\begin{proposition}\label{XiV=XiZ} Let $Z \subset \BP V$ be a non-singular variety such that
\begin{enumerate} \item[(1)] $Z$ is  linearly normal, i.e., $H^0(Z, \sO(1)) = V^*;$
\item[(2)] the variety of tangent lines of $Z$ is linearly
non-degenerate in $\BP (\wedge^2 V)$; and \item[(3)]  $H^0(Z,
ad(T(Z)) \otimes \sO(1)) =0$. \end{enumerate} Then $\Xi_V =
\Xi_Z$. \end{proposition}

\begin{proof} From the definition of $\Xi_Z$ in Proposition \ref{XiZ}, there exists a natural homomorphism
$$j: \Xi_Z \to H^0(Z, T(Z) \otimes T^*(Z) \otimes \sO(1)).$$ We claim that $j$ is
injective.
 If $\sigma \in \Xi_Z$ satisfies $j(\sigma) =0$, then
for any $u \in \hat{Z}$ and any $v \in T_u(\hat{Z})$, $\sigma(u,
v) =0$. But $\{ u \wedge v, u \in \hat{Z}, v \in T_u(\hat{Z})\}$
generates $\wedge^2 V$ by the condition (2). Thus $\sigma =0$ and
$j$ is injective.

 By the conditions (1) and (3), $$
H^0(Z, T(Z) \otimes T^*(Z) \otimes \sO(1)) =
 H^0(Z, ad(T(Z)) \otimes \sO(1)) \oplus H^0(Z, \sO(1)) \cong
V^*.$$  Note that $\Xi_V \subset \Xi_Z$. Thus the injection
$$ V^* \cong \Xi_V \subset \Xi_Z \stackrel{j}{\rightarrow} H^0(Z,
T(Z) \otimes T^*(Z) \otimes \sO(1)) \cong V^*$$ must be an
isomorphism and $\Xi_V = \Xi_Z$.
\end{proof}

From Theorem \ref{XiV},  Proposition \ref{XiZ} and Proposition
\ref{XiV=XiZ}, we have the following differential geometric
result.

\begin{theorem}\label{DG} Let $Z \subset \BP V$ be a non-singular projective
subvariety satisfying the conditions of Proposition \ref{XiV=XiZ}.
Let $\sC \subset \BP T(M)$ be a $Z$-isotrivial cone structure with
an adapted frame $\omega$. If the conic connection $\Gamma$
induced by $\omega$ on $\sC$ is a characteristic connection, then
the structure function $\sigma^{\omega}$ is conformally closed. In
particular,  the cone structure is locally flat by Proposition
\ref{flat}.
\end{theorem}

\section{Proof of Theorem \ref{Main T}}

 The following is a restatement of Proposition
3.1.2 of \cite{HM99} or Proposition 8 of \cite{HM04}.

\begin{proposition}\label{HM} Let $X$ be a projective manifold with a family of minimal rational
curves satisfying Assumption \ref{assumption}. There exists a
connected open subset $M \subset X$ such that the varieties of
minimal rational tangents define a cone structure $\sC \subset \BP
T(M)$. Then $\sC$ admits a characteristic  connection $\sF$. In
fact, the leaves of $\sF$ are given by tangent vectors to members
of $\sK$.
\end{proposition}

We remark that Assumption \ref{assumption} is crucial here. The
assumption that $\sC_x$ is non-singular implies that the tangent
map in \cite{HM04} is an embedding, which shows the existence of
$\sF$ as a line subbundle of $T(\sC)$ with locally free quotient.

 Now we can finish
the proof of Theorem \ref{Main T} as follows. From Proposition
\ref{HM}, the cone structure has a characteristic connection $\sF
\subset T(\sC)$. Choose an adapted coframe $\omega$. The induced
connection $\Gamma \subset T(\sC)$ must agree with the
characteristic connection $\sF$ by Proposition \ref{unique}. By
Theorem \ref{DG}, the cone structure is locally flat. It is easy
to see that the corresponding cone structure in Example \ref{Ex}
is also locally flat. Thus  Theorem \ref{Main T} follows from
Theorem \ref{HM99}.

\bigskip
{\bf Acknowledgment} It is a pleasure to thank Eckart Viehweg for
the reference  \cite{BR}. This work had started during my visit to
MPIM-Bonn in September-October of 2008. I would like to thank
MPIM-Bonn for the support and the hospitality.

\end{document}